\documentclass[10pt]{amsart}
\usepackage{amsmath,amsthm,amssymb}
\usepackage{graphicx}
\usepackage[margin=1.5in]{geometry}

\begin{document}


\newcommand{\PGL} {\Pj\Gl_2(\R)}                           
\newcommand{\PGLC} {\Pj\Gl_2(\Cx)}                         
\newcommand{\RP} {\R\Pj^1}                                 
\newcommand{\CP} {\Cx\Pj^1}                                
\newcommand{\C} {{\mathbb C}}                              
\newcommand{\R} {{\mathbb R}}                              
\newcommand{\Z} {{\mathbb Z}}                              
\newcommand{\Pj} {{\mathbb P}}                             
\newcommand{\Sg} {{\mathbb S}}                             

\newcommand{\suchthat} {\ \ | \ \ }
\newcommand{\ore} {\ \ {\it or} \ \ }
\newcommand{\oand} {\ \ {\it and} \ \ }

\newcommand{\oM} [1] {\ensuremath{{\mathcal M}_{0,#1}(\R)}}              
\newcommand{\M} [1] {\ensuremath{{\overline{\mathcal M}}{_{0, #1}(\R)}}}    

\newcommand{\Tubeset} {\mathfrak{T}}                     

\newcommand{\D} {\Delta}
\newcommand{\Pol} {\mathcal{P}}                           
\newcommand{\GG} {\Gamma_G}                               
\newcommand{\PG} {\Pol G}                                 

%
%

\theoremstyle{plain}
\newtheorem{thm}{Theorem}
\newtheorem{prop}[thm]{Proposition}
\newtheorem{cor}[thm]{Corollary}
\newtheorem{lem}[thm]{Lemma}
\newtheorem{conj}[thm]{Conjecture}

\theoremstyle{definition}
\newtheorem*{defn}{Definition}
\newtheorem*{exmp}{Example}

\theoremstyle{remark}
\newtheorem*{rem}{Remark}
\newtheorem*{hnote}{Historical Note}
\newtheorem*{nota}{Notation}
\newtheorem*{ack}{Acknowledgments}
\numberwithin{equation}{section}

\title {A realization of graph-associahedra}


\author{Satyan L.\ Devadoss}
\address{S.\ Devadoss: Williams College, Williamstown, MA 01267}
\email{satyan.devadoss@williams.edu}

\begin{abstract}
Given any finite graph $G$, we offer a simple realization of the graph-associahedron $\PG$ using integer coordinates.
\end{abstract}


\maketitle


\baselineskip=15pt

%
%
\section{Introduction}

Given a finite graph $G$, the graph-associahedron $\PG$ is a simple, convex polytope whose face poset is based on the connected subgraphs of $G$.  This polytope was has been studied in \cite{cd}, and has appeared in combinatorial \cite{arw, pos} and geometric contexts \cite{djs, tol}. In particular, it appears as tilings of minimal blow-ups of certain Coxeter complexes, which themselves are natural generalizations of the Deligne-Knudsen-Mumford compactification \M{n} of the real moduli space of curves \cite{dev}.

For special examples of graphs, their graph-associahedra become well-known, sometimes classical, polytopes.  For instance, when $G$ is a set of vertices, $\PG$ is the simplex.  Moreover, when $G$ is a path, a cycle, or a complete graph, $\PG$ results in the associahedron, cyclohedron, and permutohedron, respectively.
Loday \cite{lod} provided a formula for the coordinates of the vertices of the associahedron which contains the classical realization of the permutohedron.  Recently, Hohlweg and Lange \cite{hl} offer different realizations of the associahedron and cyclohedron.   We offer a realization of graph-associahedra for any graph.

%
%
\section{Convex Hull}
\subsection{}

We begin with definitions; the reader is encouraged to see \cite[Section 1]{cd} for details.

\begin{defn}
Let $G$ be a finite graph.  A \emph{tube} is a proper nonempty set of nodes of $G$ whose induced graph is a proper, connected subgraph of $G$.  There are three ways that two tubes $u_1$ and $u_2$ may interact on the graph.
\begin{enumerate}
\item Tubes are \emph{nested} if  $u_1 \subset u_2$.
\item Tubes \emph{intersect} if $u_1 \cap u_2 \neq \emptyset$ and $u_1 \not\subset u_2$ and $u_2 \not\subset u_1$.
\item Tubes are \emph{adjacent} if $u_1 \cap u_2 = \emptyset$ and $u_1 \cup u_2$ is a tube in $G$.
\end{enumerate}
Tubes are \emph{compatible} if they do not intersect and they are not adjacent.  A \emph{tubing} $U$ of $G$ is a set of tubes of $G$ such that every pair of tubes in $U$ is compatible.  A \emph{$k$-tubing} is a tubing with $k$ tubes.
\end{defn}

\begin{figure}[h]
\includegraphics[width=\linewidth]{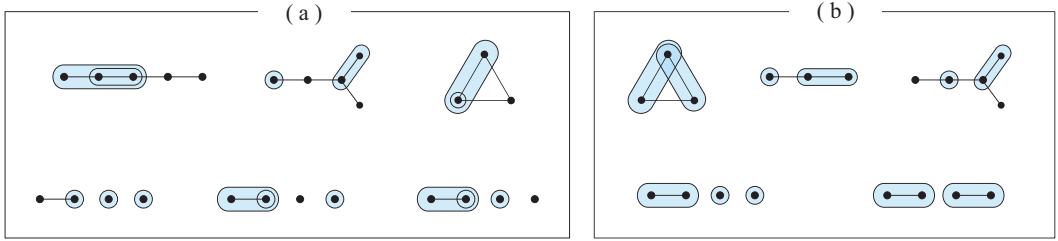}
\caption{(a) Valid tubings and (b) invalid tubings.}
\label{f:tubings}
\end{figure}

\begin{rem}
When $G$ is a disconnected graph with connected components $G_1$, \ldots, $G_k$, we place an additional restriction.  Let $u_i$ be the tube of $G$ whose induced graph is $G_i$.  Then any tubing of $G$ cannot contain all of the tubes $\{u_1, \ldots, u_k\}$.  Thus, for a graph $G$ with $n$ nodes, a tubing of $G$ can at most contain $n-1$ tubes.  Figure~\ref{f:tubings} shows examples of (a) valid tubings and (b) invalid tubings.
\end{rem}

\begin{defn}
For a graph $G$, the \emph{graph-associahedron} $\PG$ is a simple, convex polytope whose face poset is isomorphic to set of tubings of $G$, ordered such that $U \prec U'$ if $U$ is obtained from $U'$ by adding tubes.
\label{d:pg}
\end{defn}

\subsection{}

\begin{figure}[h]
\includegraphics[width=\linewidth]{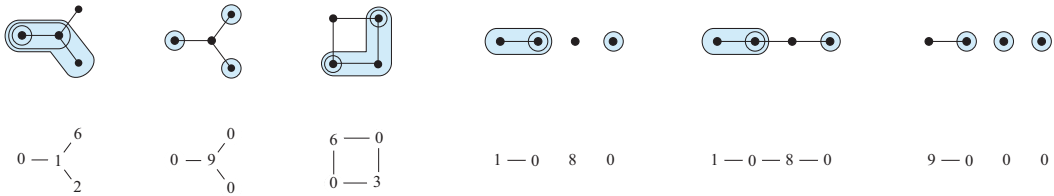}
\caption{Integer values of nodes associated to tubings.}
\label{f:coord}
\end{figure}

Let $G$ be a graph with $n$ nodes and let $M_G$ be the collection of maximal $(n-1)$-tubings of $G$. For each such tubing $U$ in $M_G$, define a map $f_U$ from the nodes of $G$ to the integers as follows:  If a node $v$ of $G$ is a tube of $U$, then $f_U(v) = 0$.  Otherwise, let $t(v)$ be the smallest tube containing $v$, 
and let all other nodes of $G$ satisfy the recursive condition
\begin{equation}
\sum_{x \in t(v)} f_U(x) = 3^{|t(v)| - 2}.
\label{e:integer}
\end{equation}
Figure~\ref{f:coord} gives some examples of integer values of nodes associated to tubings.

\noindent Let $G$ be a graph with an ordering $v_1, v_2, \ldots, v_n$ of its nodes.  Define $c: M_G \rightarrow \R^n$ where
$$c(U) = (f_U(v_1), f_U(v_2), \ldots, f_U(v_n)).$$

\begin{thm}
If $G$ is a graph with $n$ nodes, the convex hull of the points $c(M_G)$ in $\R^n$ yields the graph-associahedron $\PG$.
\label{t:hull}
\end{thm}

\noindent The proof of this is given at the end of the paper.

%
%
\section{Examples}
\subsection{Simplex}

\begin{figure}[h]
\includegraphics[width=\linewidth]{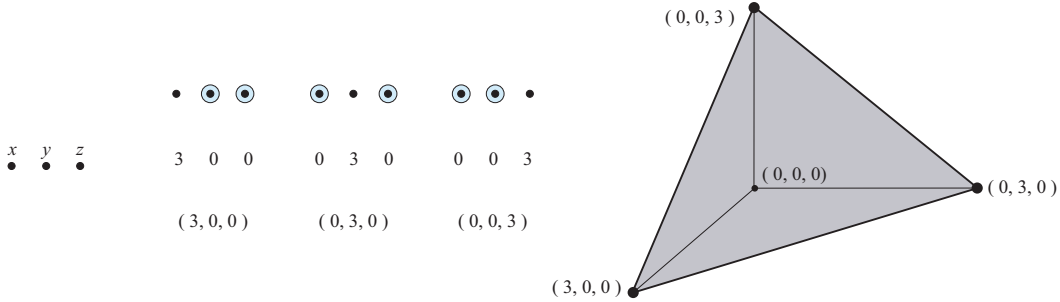}
\caption{The maximal tubings of $G$ and its convex hull, resulting in the simplex.}
\label{f:simplex}
\end{figure}

Let $G$ be the graph with $n$ (disjoint) nodes.  The set $M_G$ of maximal tubings has $n$ elements, each corresponding to choosing $n-1$ out of the $n$ possible nodes.  An element of $M_G$ will be assigned a point in $\R^n$ consisting of zeros for all coordinates except one with value $3^{n-2}$.  Due to Theorem~\ref{t:hull}, $\PG$ is the convex hull of the $n$ vertices in $\R^n$ yielding the $(n-1)$-simplex.  Figure~\ref{f:simplex} shows this when $n = 3$, resulting in the $2$-simplex in $\R^3$.

\subsection{Permutohedron}

\begin{figure}[h]
\includegraphics[width=\linewidth]{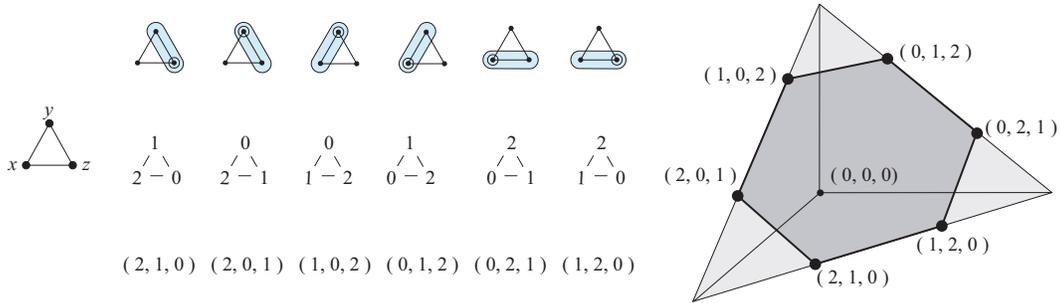}
\caption{The maximal tubings of $G$ and its convex hull, resulting in the permutohedron.}
\label{f:permutohedron}
\end{figure}

Let $G$ be the complete graph on $n$ nodes.  Each maximal tubing of $G$ can be seen as a sequential nesting of all $n$ nodes.  In other words, they are in bijection with permutations on $n$ letters.  The elements of $M_G$ will be assigned coordinate values based on all permutations of $\{0, 1, \ldots, 3^{n-2}-3^{n-3}\}$.  Theorem~\ref{t:hull} shows $\PG$ as the convex hull of the $n!$ vertices in $\R^n$, resulting in the permutohedron. Figure~\ref{f:permutohedron} shows this when $n = 3$, resulting in the hexagon, the two-dimensional permutohedron.

\subsection{Associahedron}

\begin{figure}[h]
\includegraphics[width=\linewidth]{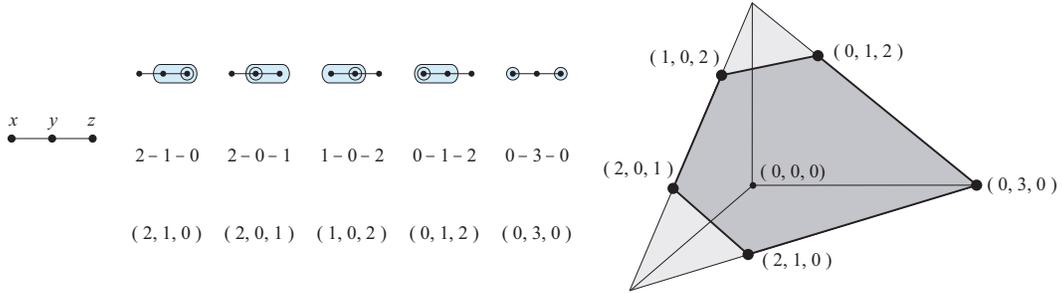}
\caption{The maximal tubings of $G$ and its convex hull, resulting in the associahedron.}
\label{f:associahedron}
\end{figure}

Let $G$ be an $n$-path.   The number of such maximal tubings is in bijection with the Catalan number $c_n$.   Due to Theorem~\ref{t:hull}, the convex hull of these vertices in $\R^n$ yields the $(n-1)$ dimension associahedron.  Stasheff originally defined the associahedron for use in homotopy theory in connection with associativity properties of $H$-spaces \cite{sta}. Figure~\ref{f:associahedron} shows this when $n = 3$, resulting in the pentagon, the two-dimensional associahedron.

\subsection{Cyclohedron}

Let $G$ be an $n$-cycle. In this case, the number of maximal tubings is the type B Catalan number $\binom{2n-2}{n-1}$.   Theorem~\ref{t:hull} shows $\PG$ as the cyclohedron, a polytope originally manifested in the work of Bott and Taubes in relation to knot and link invariants \cite{bt}.  Figure~\ref{f:permutohedron} shows this when $n = 3$, since both the permutohedron and cyclohedron are identical in dimension two.

%
%
\section{Constructing the Graph-Associahedron}
\subsection{}

For a graph $G$ with $n$ nodes $v_1, \ldots v_n$, let $\D$ be the $(n-1)$-simplex in which each facet (codimension $1$ face) corresponds to a particular node of $G$.  Thus, each proper subset of nodes of $G$ corresponds to a unique face of $\D$, defined by the intersection of the faces associated to those nodes.  The following construction of the graph-associahedron is based on truncations of a simplex.

\begin{thm}\cite[Section 2]{cd}
For a given graph $G$, truncating faces of $\D$ which correspond to $1$-tubings in increasing order of dimension results in $\PG$.
\label{t:pg}
\end{thm}

\begin{figure}[h]
\includegraphics[width=.95\linewidth]{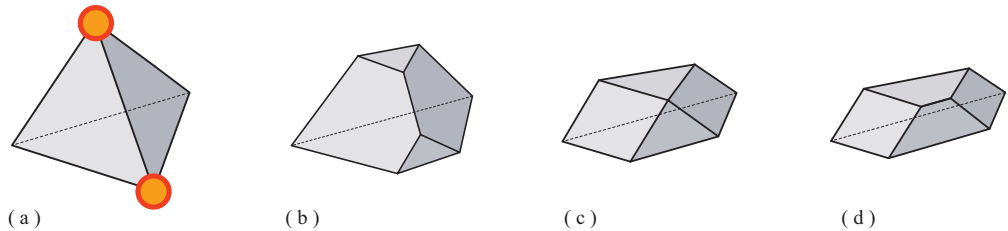}
\caption{Iterated truncations of the $3$-simplex based on an underlying graph.}
\label{f:deepcuts}
\end{figure}

\noindent   Indeed, truncations should not only be in increasing order of dimension (certain vertices of $\D$ are truncated first, and then the edges, and so forth), but they should also not form ``deep cuts''.  Consider Figure~\ref{f:deepcuts} as an example.  Part (a) shows a $3$-simplex with two vertices marked for truncation;  part (b) shows appropriate truncations of the vertices, with (c) and (d) showing inappropriate cuts which are too deep.

\begin{rem}
In order to recover Loday's elegant construction of the classical permutohedron as part of the associahedron, we simply use the following recursive definition of $f_U$:
$$\sum_{x \in t(v)} f_U(x) = \binom{|t(v)|+1}{2}.$$
Although this works for the associahedron, it fails for graph-associahedra in general.  The reason for this is that the cuts needed to construct the polytopes are too deep.
\end{rem}

\begin{figure}[h]
\includegraphics[width=.95\linewidth]{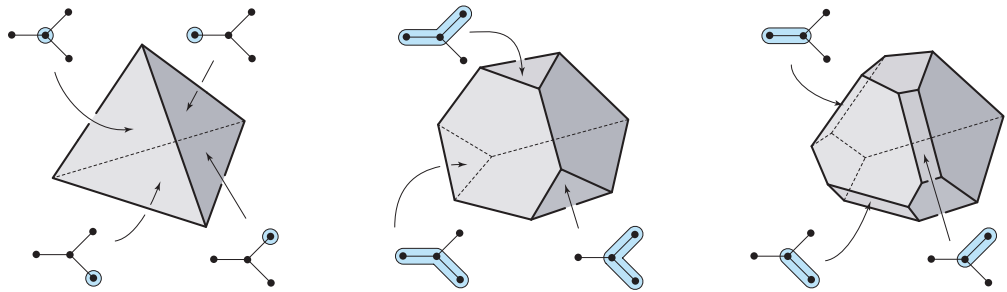}
\caption{Iterated truncations of the $3$-simplex based on an underlying graph.}
\label{f:d4}
\end{figure}

Figure \ref{f:d4} shows a tetrahedron truncated according to a graph, resulting in $\PG$.  Note that its facets are labeled with $1$-tubings.  One can verify that the edges correspond to all possible $2$-tubings and the vertices to $3$-tubings.

\subsection{}

We are now in position to prove Theorem~\ref{t:hull}.  This is influenced by the work of Stasheff and Schnider \cite[Appendix B]{sta2}.  

\begin{proof}[Proof of Theorem~\ref{t:hull}]
Consider the affine hyperplane $H$ defined by
\begin{equation}
\sum x_i = 3^{n-2}.
\label{e:hyp}
\end{equation}
The intersection of the quadrant $\{(x_1, \ldots, x_n) \ | \ x_i \geq 0\}$ with $H$ yields a standard $(n-1)$-simplex $\D$.  Let $G_u$ be the set of all $1$-tubings of $G$, where $G_u^i$ be the set of $1$-tubings containing $i$ nodes.
The faces of $\D$ which need to truncated correspond to the $1$-tubings $G_u^i$, where $i \geq 2$.  Let $u = \{v_{i_1}, \ldots, v_{i_k}\}$ be a $1$-tubing in $G_u^k$; note that this corresponds to a $n-1-k$ face of $\D$, seen as the intersection of the hyperplane
$$\sum_{v_i \in u} x_i = 0$$
of $\R^n$ with $\D$.  Truncate this face with the hyperplane
\begin{equation}
\sum_{v_i \in u} x_i = 3^{k-2}.
\label{e:trunc}
\end{equation}
We claim that this collection of hyperplanes, one for each element of $G_u^i$, results in $\PG$.

By Theorem~\ref{t:pg} above, the appropriate faces of $\D$ have been truncated, one for each $1$-tubing.  However, we need to show the any two cuts of a given dimension are not deep;  that is, their corresponding hyperplanes must not intersect in $H$.  This is done by induction.  Two vertices of $\D$ which are truncated correspond to $1$-tubings in $G_u^{n-1}$, say $u = \{v_1, \ldots, v_{n-2}, v_{n-1}\}$ and $u' = \{v_1, \ldots, v_{n-2}, v_n\}$.  These hyperplanes cannot intersect in $H$ since Eqs.~\eqref{e:hyp} and \eqref{e:trunc} show
$$\sum x_i \ = \ 3^{n-2} \ > \ 3^{n-3} + 3^{n-3} \ = \ \sum_{v_i \in u} x_i + \sum_{v_i \in u'} x_i.$$
In general, let $u_*$ be in $G_u^{n-1-k}$, a $k$-dimensional face of $\D$ that is truncated.  Let $u$ and $u'$ be two $(k+1)$-dimensional faces of $\D$ in $G_u^{n-2-k}$ which are incident to $u_*$.  The cuts $u$ and $u'$ will not be deep with respect to $u_*$.  To see this, notice that the nodes of $u$ and $u'$ are contained in $u_*$.  Thus, in $H$ the hyperplanes of $u$ and $u'$ cannot intersect in $u_*$ since
$$\sum_{v_i \in u_*} x_i \ = \ 3^{n-3-k} \ > \ 3^{n-2-k} + 3^{n-2-k} \ = \ \sum_{v_i \in u} x_i + \sum_{v_i \in u'} x_i.$$

Recall that each vertex of $\PG$ corresponds to a $(n-1)$-tubing $T$ of $G$.  This, in turn, corresponds to the intersection of the $n-1$ hyperplanes of \eqref{e:trunc} for each $1$-tubing of $T$.  In particular, a  tube containing one node assigns the value $0$ to that node; these are incident to the original facets of $\D$.  Thus Eq.~\eqref{e:integer} is satisfied inductively. 
\end{proof}


\begin{ack}
We thank the mathematics department at the Ohio State university for the opportunity to visit during the 2005-2006 academic year, where this work was completed.  We are also grateful to Igor Pak, for motivating conversations, and to the NSF for partially supporting this work with grant DMS-0310354.
\end{ack}

%
%
\bibliographystyle{amsplain}

\begin{thebibliography}{XXX}
\baselineskip=12pt

\bibitem[1]{arw} F.\ Ardila, V.\ Reiner, L.\ Williams. Bergman complexes, Coxeter arrangements, and graph associahedra, \emph{Seminaire Lotharingien de Combinatoire} {\bf 54A}(2006).

\bibitem[2]{bt} R.\ Bott and C.\ Taubes. On the self-linking of knots, \emph{J.\ Math.\ Phys.} {\bf 35} (1994) 5247-5287.

\bibitem[3]{cd} M.\ Carr and S.\ L.\ Devadoss.  Coxeter Complexes and graph-associahedra, \emph{Topology and its Appl.} {\bf 153} (2006) 2155-2168.

\bibitem[4]{djs} M.\ Davis, T.\ Januszkiewicz, R.\ Scott.  Fundamental groups of blow-ups, \emph{Advances in Math.} {\bf 177} (2003) 115-179.

\bibitem[5]{dev} S.\ Devadoss. Tessellations of moduli spaces and the mosaic operad,  {\em Contemp.\ Math.} {\bf 239} (1999) 91-114.

\bibitem[6]{hl} C.\ Hohlweg and C.\ Lange.  Realizations of the associahedron and cyclohedron, preprint math.CO/0510614.

\bibitem[7]{lod} J.-L.\ Loday.  Realization of the Stasheff polytope, \emph{Archiv der Mathematik} {\bf 83} (2004) 267-278.

\bibitem[8]{pos} A.\ Postnikov.  Permutohedra, associahedra, and beyond, preprint math.CO/0601339.

\bibitem[9]{sta} J.\ D.\ Stasheff. Homotopy associativity of $H$-spaces, \emph{Trans.\ Amer.\ Math.\ Soc.} {\bf 108} (1963) 275-292.

\bibitem[10]{sta2} J.\ D.\ Stasheff (Appendix B coauthored with S.\ Shnider). From operads to ``physically'' inspired theories, \emph{Contemp.\ Math.} {\bf 202} (1997) 53-81.

\bibitem[11]{tol} V.\ Toledano-Laredo.  Quasi-Coxeter algebras, Dynkin diagram cohomology and quantum Weyl groups, preprint math.QA/0506529.

\end{thebibliography}

\end{document}